\documentclass[12pt]{amsart}
\usepackage{amssymb}

 \rm

\renewcommand{\(}{\left(}
\renewcommand{\)}{\right)}

\newcommand{\x}{\times}
\renewcommand{\bar}{\overline}
\newcommand{\abs}[1]{\left\lvert#1\right\rvert}
\newcommand{\norm}[1]{\left\lVert#1\right\rVert}
\newcommand{\st}{\:|\:}

\renewcommand{\phi}{\varphi}

\renewcommand{\Re}{{\mathrm{Re}\,}}

\renewcommand{\span}{\mathrm{span}}

\newcommand{\CC}{\mathbb{C}}
\newcommand{\RR}{\mathbb{R}}
\newcommand{\ZZ}{\mathbb{Z}}
\newcommand{\TT}{\mathbb{T}}
\newcommand{\DD}{\mathbb{D}}

\renewcommand{\H}{\mathcal{H}}
\newcommand{\BH}{\mathcal{B}(\H)}
\newcommand{\Mob}{\textrm{M\"ob}}
\newcommand{\M}{\mathcal{M}}
\newcommand{\A}{\mathcal{A}}

\theoremstyle{plain}
\newtheorem{thm}{Theorem}[section]
\newtheorem{lem}[thm]{Lemma}
\newtheorem{prop}[thm]{Proposition}

\theoremstyle{definition}

\theoremstyle{remark}
\newtheorem{rem}[thm]{Remark}

\title{Inductive algebras and homogeneous shifts}

\author{Amritanshu Prasad}
\address{The Institute of Mathematical Sciences, CIT campus Taramani,
Chennai 600 113.}

\author{M.~K.~Vemuri}
\address{Chennai Mathematical Institute, Plot H1, SIPCOT IT Park,
Padur~PO, Siruseri 603 103.}

\begin{document}

\begin{abstract}
Inductive algebras for the irreducible unitary representations of the
universal cover of the group of unimodular two by two matrices are
classified.  The classification of homogeneous shift operators is
obtained as a direct consequence.  This gives a new approach to the
results of Bagchi and Misra.
\end{abstract}

\maketitle


\section{Introduction}\label{S:intro}

Let $G$ be a separable locally compact group and let $R$ be a
strongly continuous representation of $G$ on a separable Hilbert space
$\H$.  Let $\BH$ denote the algebra of bounded operators on $\H$.  An
{\em inductive algebra} is a strong-operator closed abelian
sub-algebra of $\BH$ that is normalized by $R(G)$.  If we wish to
emphasize the dependence on $R$, we use the term $R$-inductive algebra.
The inductive algebras for the irreducible unitary representations of
$\mathrm{SL}(2,\RR)$ were classified in \cite{sl2r}.  In this work, we
classify the inductive algebras for the irreducible unitary representations
of the universal cover of $\mathrm{SL}(2,\RR)$.  We use this result
to give another proof of the classification of the homogeneous shifts
originally obtained by Bagchi and Misra \cite{Bagchi-Misra}.

Let $\DD$ and $\TT$ denote the unit disc and the unit circle in the
complex plane $\CC$ respectively.  Let $\Mob$ denote the M\"obius
group (the group of biholomorphic automorphisms of $\DD$).  Thus
$\Mob = \{\phi_{\alpha,\beta} \st \alpha \in \TT, \beta \in \DD \}$,
where
$$
\phi_{\alpha,\beta}(z) = \alpha \frac{z-\beta}{1-\bar{\beta}z},
\qquad z \in \DD.
$$
The correspondence $(\alpha, \beta) \mapsto \phi_{\alpha,\beta}$ is
bijective and identifies $\Mob$ with the manifold $\TT \x \DD$.
With this differential structure, $\Mob$ is a Lie group isomorphic
to $\mathrm{PSL}(2, \RR)$.  An operator $T \in \BH$ is called
{\em homogeneous} if $\phi(T)$ is unitarily equivalent to $T$ for those
$\phi \in \Mob$ which are holomorphic on the
spectrum of $T$.  The homogeneous weighted shift operators
(the homogeneous shifts, for brevity) were classified in
\cite{Bagchi-Misra}.

It was noticed by Gadadhar Misra that the operators which
(topologically) generate the inductive algebras found in \cite{sl2r}
are weighted shifts whose weight-sequences were the same as
{\em some} of the operators found in \cite{Bagchi-Misra}.  As explained
in \cite{Bagchi-Misra}, every homogeneous shift is (non-uniquely)
associated to a {\em simple} representation of the universal cover
$\widetilde{\Mob}$ of
$\Mob$.  Furthermore, it turns out that the closure in the strong-operator
topology of the sub-algebra generated by the set
$\{\phi(T) \st \phi \in \Mob\}$ is inductive for this representation
(see Prop. \ref{P:homogeneous-inductive}).  Thus a classification of the
inductive algebras for the simple representations of $\widetilde{\Mob}$
would yield {\em all} the operators that were found in \cite{Bagchi-Misra},
and in fact an independent proof of the classification (Theorem 5.2 in
\cite{Bagchi-Misra}).

Here, we find the inductive algebras only for the irreducible
unitary representations of $\widetilde{\Mob}$.  That still leaves
the reducible simple representations, i.e. those of the form
$D_{2-\lambda}^- \oplus D_\lambda^+$, $0 < \lambda < 2$.  The
{\em assumption} that $T$ is a (scalar) shift allows us to
patch together information contained in the inductive algebras
for $D_{2-\lambda}^-$ and $D_\lambda^+$ to determine the
homogeneous operators associated to
$D_{2-\lambda}^- \oplus D_\lambda^+$.

Every homogeneous operator is a block shift \cite{Bagchi-Misra},
but the problem of classifying homogeneous operators remains
open.  This work addresses the case where the blocks are one
dimensional.  Extending it elegantly to the case of higher dimensional
blocks would require the classification of inductive algebras for
reducible representations
(the trick for treating $D_{2-\lambda}^- \oplus D_\lambda^+$ does
not work in general).
Classifying inductive algebras for reducible representations appears
to be a much harder problem (akin to classifying all abelian sub-algebras
of a matrix algebra).

The homogeneous operators belonging to a Cowen-Douglas class were
recently classified (see \cite{Wilkins}, \cite{Koranyi-Misra-1},
\cite{Koranyi-Misra-2}).

In section \ref{S:UD} we recall the relevant notation from
\cite{Bagchi-Misra} and \cite{sl2r} and describe the irreducible
representations of $\widetilde{\Mob}$ in a convenient way.
In section \ref{S:IA}, we determine the inductive algebras for each of
these representations.  In section \ref{S:HS}, we use the results of
section \ref{S:IA} to classify the homogeneous shifts.

\section{The unitary dual of $\widetilde{\Mob}$}\label{S:UD}

Henceforth let $G=\widetilde{\Mob}$ and
let $\pi:G \to \Mob$ denote the universal covering map.  If $g \in G$
let $\phi_g$ denote $\pi(g)$ thought of as a function $\DD \to \DD$.
Since $G$ is connected and simply-connected, for each $\eta \in \CC$
there is a unique smooth branch of the function
$(\phi_{g^{-1}}'(z))^{\eta}$ such that $(\phi_1'(z))^{\eta} = 1$.

Recall that if $\phi \in \Mob$ and $\phi^*(z) := \bar{\phi(\bar{z})}$
then $\phi \mapsto \phi^*$ is an automorphism of $\Mob$ (see
equation (2.1) in \cite{Bagchi-Misra}).  We let $g \mapsto g^*$ be the
lift of this automorphism to $G$.  Then
$$
\phi_{g^*} = \phi_g^*.
$$
If $R$ is a representation of $G$, then $R^\#$ is the
representation defined by
$$
R^\#(g) = R(g^*), \qquad g \in G.
$$
This notation is an extension of equation (2.4) in \cite{Bagchi-Misra}.

Let $\lambda \in \RR$, $\mu \in \CC$ and $I \in \{ \ZZ, \ZZ^+ \}$, and
assume $\mu=0$ if $I = \ZZ^+$.  Let $\M(\ZZ)$ (resp.\ $\M(\ZZ^+)$)
denote the functions which have a holomorphic extension to some
neighborhood of $\TT$ (resp. $\bar{\DD}$).

For $g \in G$, define
$R(g)=R_{\lambda,\mu}(g):\mathcal{M}(I) \to \mathcal{M}(I)$
by
$$
\begin{aligned}
(R(g)F)(z)
:= & \; (\phi_{g^{-1}}'(z))^{\lambda/2} \abs{\phi_{g^{-1}}'(z)}^\mu
        F(\phi_{g^{-1}}(z)) \\
=  & \; (\phi_{g^{-1}}'(z))^{(\lambda + \mu)/2}
        \bar{(\phi_{g^{-1}}'(z))^{\mu/2}}
        F(\phi_{g^{-1}}(z)).
\end{aligned}
$$
Then $g \mapsto R(g)$ is a representation of $G$ on $\mathcal{M}(I)$.

Let $f_n(z) = z^n$.  If we give $\M(I)$ the $C^\infty$ (Frechet)
topology, then the linear span of $\{ f_n \}_{n \in I}$ is dense in
$\M(I)$.  Let
\begin{equation}\label{E:norms}
\norm{f_n}^2 = \frac{\Gamma(1-\mu+n)}{\Gamma(\lambda+\bar{\mu}+n)},
\qquad n \in I,
\end{equation}
when $\lambda$ and $\mu$ are such that the expression on the right
is real and positive.  Extend
$\norm{\cdot}$ to $\span\{ f_n \}_{n \in I}$ as a norm.  Since
$\norm{f_n}$ grows at most like a polynomial in $\abs{n}$ as
$\abs{n} \to \infty$, it follows that $\norm{\cdot}$ is uniformly
continuous, and thus extends to $\M(I)$.  Let $\H = \H^{\lambda,\mu}$
be the completion of $\M(I)$ under this norm.  Then $R$ extends to
a unitary representation of $G$ on $\H$.
Puk{\'a}nszky \cite{Pukanszky} showed that every irreducible unitary
representation of $G$ is unitarily equivalent either to a representation
of this type or to the composition of one with the $*$-automorphism.
We recall his taxonomy in terms of our parameters $(I, \lambda, \mu)$:

\begin{description}
\item[Holomorphic discrete series] $D_\lambda^+ = R_{\lambda, 0}$ where
$I=\ZZ^+$ and $\lambda>0$.

\item[Anti-holomorphic discrete series] $D_\lambda^- = (D_\lambda^+)^{\#}$,
$\lambda > 0$.

\item[Principal series] $R_{\lambda, \mu} $, $I=\ZZ$, $\lambda \in (-1, 1]$ and
$\Re \mu = \frac{1-\lambda}{2}$.

\item[Complementary series] $R_{\lambda, \mu}$, $I=\ZZ$, $\lambda \in (-1,1)$
and $\mu \in (0,1) \cap (-\lambda, 1-\lambda)$.

\end{description}

Within the principal and complementary series, there are unitary
equivalences between $R_{\lambda, \mu}$ and $R_{\lambda, 1-\lambda-\mu}$
and $R_{1,0} \cong D_1^+ \oplus D_1^-$.  Otherwise, these representations
are irreducible, inequivalent and cover the unitary dual of
$G$.  We remark that the usage ``Discrete series'' here is a historical
accident, and in fact $D_\lambda^\pm$ do not embed in $L^2(G)$.

\section{The inductive algebras}\label{S:IA}

Let $R$ be one of the irreducible unitary representations $R_{\lambda,\mu}$
introduced in section \ref{S:UD}.
For $g\in G$, define $\tilde\kappa(g):\BH\to\BH$ by
$T \mapsto R(g)TR(g)^{-1}$.  If $g\in \ker(\pi)$ then
$R(g)$ is a scalar, so $\tilde\kappa(g)$ is trivial.  So
$\tilde\kappa$ descends to a representation $\kappa$ of $\textrm{M\"ob}$.
The remarks about $\kappa$ at the beginning of section 3 in
\cite{sl2r} continue to be valid.

Let $K=\{ \phi_{\alpha,0} \st \alpha \in \TT \} \subseteq \textrm{M\"ob}$.
Then $K$ is a maximal compact subgroup.  For $l \in \ZZ$, let
$\chi_l: K \to \CC^*$ be defined by $\chi_l(\phi_{\alpha, 0}) = \alpha^l$.
Then $\widehat{K}=\{\chi_l\}_{l\in\ZZ}$.

Let $\mathfrak{g}$ denote the Lie algebra of $\Mob$.  We let $\exp$
denote the exponential map from $\mathfrak{g}$ to any of its associated
Lie groups.  The precise group will be clear from context. Let
$h, L, M \in \mathfrak{g}$ be such that
$$
\begin{aligned}
\exp{th} &= \phi_{e^{2it},0}, \\
\exp{tL} &= \phi_{1, -\tanh t}, \\
\exp{tM} &= \phi_{1, -i\tanh t},
\end{aligned}
$$
and $e=\frac12 (L-iM)$ and $f=\frac12 (L+iM)$.  Then
$h, e, f \in \mathfrak{g} \otimes \CC \cong \mathfrak{sl}(2, \CC)$.

Observe that
$$
\{ F \in \H \st R(\exp(th))F = e^{-i(2n+\lambda)t} F \} = \CC f_n,
\qquad n \in I.
$$
Differentiating along one-parameter subgroups (and then forming
complex linear combinations) gives
$$
\begin{aligned}
R(e)f_n
= & \;
\begin{cases}
(\mu - n) \, f_{n-1}, & n \in I \cap (I+1), \\
0                     & \text{otherwise,}
\end{cases} \\
R(f)f_n
= & \; (\lambda + \mu + n) \, f_{n+1}, \qquad n \in I.
\end{aligned}
$$

Let $\A \subseteq B(\mathcal{H})$ be an $R$-inductive
algebra.    For $m \in \ZZ$, define the vector spaces
$$
\A_m = \{ T \in \A \st \kappa(g)T = \chi_m(g) T, \, g \in K \}.
$$
Then by the Peter-Weyl theorem
$$
\bigoplus_{m \in \ZZ} \A_m
$$
is sequentially dense in $\A$ (see section 3 of \cite{sl2r}).

If $T \in \A_m$ then
$$
\begin{aligned}
R(\exp(th)) T f_n
= & \; \kappa(\exp(th))T R(\exp(th)) f_n \\
= & \; \chi_m(\exp(th)) T R(\exp(th)) f_n \\
= & \; e^{i(2m - 2n - \lambda)t} T f_n,
\end{aligned}
$$
so $T f_n \in \CC f_{n-m}$.  So 
$$
T f_n =
\begin{cases}
a_n f_{n-m} & n \in I \cap (I+m), \\
0           & \text{otherwise},
\end{cases}
$$
for some $a_n \in \CC$.  A similar
calculation shows that conversely if $T$ has this form, then $T \in \A_m$.

Let $\H^\infty$ and $\A^\infty$ denote the smooth vectors in $\H$ and
$\A$ respectively.  Recall that $v$ is smooth if the orbit map
$g \mapsto g \cdot v$ is a smooth mapping on $G$.  We recall the
following facts from section 3 of \cite{sl2r}, and use them without
mention in the calculations that follow.

\begin{itemize}

\item $\A^\infty$ is invariant under $\kappa(G)$ as well as
$\kappa(\mathfrak{g})$.

\item If $T \in \A^\infty$ and $F \in \H^\infty$ then $TF \in \H^\infty$.

\item $\A_m \cap \A^\infty$ is sequentially dense in $\A_m$.

\end{itemize}

For $T \in \A^\infty$, define $T_e = [R(e), T]$ and $T_f = [R(f), T]$.
Then $T_e, T_f \in \A^\infty$.  If $T \in \A_m \cap \A^\infty$, then
$T_e \in \A_{m+1}$ and $T_f \in \A_{m-1}$.  Moreover,

\begin{equation}\label{E:TeTf}
\begin{aligned}
T_e f_n
= & \; ((\mu - n + m) a_n - (\mu - n) a_{n-1}) f_{n-m-1},
\qquad n \in (I+1) \cap (I+m+1), \\
T_f f_n
= & \; ((\lambda + \mu + n - m) a_n - (\lambda + \mu + n) a_{n+1}) f_{n-m+1},
\qquad n \in I \cap (I+m).
\end{aligned}
\end{equation}
For other values of $n$ these formulas degenerate, and we will have to
consider cases, when the need arises.

\begin{lem} \label{L:strong-Schur}
$\A_0=\CC \mathrm{I}$
\end{lem}

\begin{proof}
It is clear that $\CC \mathrm{I} \subseteq \A_0$.  Let
$T \in \A_0 \cap \A^\infty$.  Then
$$
\begin{aligned}
0
= & \; [T, T_e] f_n \\
= & \; -(\mu - n)(a_{n-1} - a_n)^2 f_{n-1}
\quad\text{for all $n \in I+1$}.
\end{aligned}
$$
Since $\mu \not\in I+1$ for any of the representations under consideration,
it follows that $a_{n+1}=a_n$ for all $n \in I$.  So
$\A_0 \cap \A^\infty = \CC I$, which is dense, and being finite dimensional,
closed in $\A_0$.
\end{proof}

\begin{lem}\label{L:up}
For $m \ne 0$, $\A_{m-1}=0 \implies \A_m=0.$
\end{lem}

\begin{proof}
Suppose $m \ne 0$, $\A_{m-1}=0$ and $T\in \A_m$.  Write
$$
T f_n = a_n f_{n-m}, \qquad n \in I \cap (I+m).
$$
Since $T_f \in \A_{m-1}$, we have

\begin{equation}\label{E:rr2}
(\lambda + \mu + n - m) a_n - (\lambda + \mu + n) a_{n+1} = 0,
\qquad n \in I \cap (I+m).
\end{equation}

Since $-(\lambda+\mu) \not\in I$ for any of the representations under
consideration, it follows that if $a_n = 0$ for some $n$ then $a_n=0$
for all $n$.

Since $T$ commutes with $T_e$, we have

$$
\begin{aligned}
0
=& \; [T,T_e]f_n \\
= & \; (-(\mu - n) a_{n-1} a_{n-m-1} + 
       2 (\mu - n + m) a_n a_{n-m-1} -
       (\mu - n + 2 m) a_n a_{n-m})
       f_{n-2m-1} \\
= & \; \frac{\(
             \begin{gathered}
             -(\mu - n)(\lambda + \mu + n - 1) \\
             + 2 (\mu - n + m)(\lambda + \mu + n - m - 1) \\
             - (\mu - n + 2m)(\lambda + \mu + n - 2m - 1)
             \end{gathered}
             \)}
            {(\lambda + \mu + n - m - 1)}
            a_n a_{n-m-1}
       f_{n-2m-1},
       \qquad\text{(for large $n$, by (\ref{E:rr2}))}, \\
= & \; \frac{2m^2}{(\lambda + \mu + n - m - 1)}
            a_n a_{n-m-1}
       f_{n-2m-1}.
\end{aligned}
$$
So for large $n$ either $a_n=0$ or $a_{n-m-1}=0$.  Now (\ref{E:rr2})
implies $a_n=0$ for all $n \in I \cap (I+m)$.
\end{proof}

\begin{lem}\label{L:down}
For $m \ne 0$, $\A_{m+1}=0 \implies \A_m=0$.
\end{lem}

\begin{proof}

Suppose $m \ne 0$, $\A_{m+1}=0$ and $T\in \A_m$.  Write
$$
T f_n = a_n f_{n-m}, \qquad n \in I \cap (I+m).
$$

Since $T_e \in \A_{m+1}$, we have

\begin{equation}\label{E:rr1}
(\mu - n + m) a_n - (\mu - n) a_{n-1} = 0, \qquad n \in (I+1) \cap (I+m+1).
\end{equation}
Since $\mu \not\in \ZZ$ for any of the representations under consideration,
it follows that if $a_n = 0$ for some $n$ then $a_n=0$ for all $n$.

Since $T$ commutes with $T_f$, we have
$$
\begin{aligned}
0
=& \; [T,T_f]f_n \\
=& \; (-(\lambda + \mu + n)a_{n+1} a_{n-m+1} +
      2 (\lambda + \mu + n -m) a_n a_{n-m+1} -
      (\lambda + \mu + n - 2m)a_n a_{n-m}) f_{n-2m+1} \\
=& \; \frac{\(
            \begin{gathered}
            -(\lambda + \mu + n)(\mu - n - 1) \\
            +2 (\lambda + \mu + n - m)(\mu - n + m - 1) \\
            -(\lambda + \mu + n - 2m)(\mu - n + 2m - 1)
            \end{gathered}
            \)}
           {\mu - n + m - 1}
      a_n a_{n-m+1} f_{n-2m+1},
      \qquad\text{(for large $n$ by (\ref{E:rr1}))}, \\
=& \; \frac{2m^2 a_n a_{n-m+1}}{\mu - n + m - 1} f_{n-2m+1},
\end{aligned}
$$
So for large $n$, either $a_n=0$ or $a_{n-m+1}=0$.  Now (\ref{E:rr1})
implies $a_n=0$ for all $n \in I \cap (I+m)$.
\end{proof}

One possibility is that $\A = \CC \mathrm{I}$.  Henceforth let us assume
that $\A \ne \CC \mathrm{I}$.  Then by the Peter-Weyl theorem there
exists $m \ne 0$ such that $\A_m \ne 0$.  By Lemmas \ref{L:up} and
\ref{L:down}, it follows that either $\A_1 \ne 0$ or $\A_{-1} \ne 0$.

To analyze the situation further, it is now best to consider cases.

\smallskip
\noindent{\bf Case $I=\ZZ^+$ ($\mu = 0$):}

Suppose $T \in \A_{-1} \cap \A^\infty$.  Then $T f_n = a_n f_{n+1}$ for
some $a_n \in \CC$, $n=0,1,2, \dots$.  Since $T_e \in \A_0$, by
Lemma \ref{L:strong-Schur} there exists $b \in \CC$ such that
$T_e = -b \mathrm{I}$.  We compute 
$$
T_e f_0 = R(e) T f_0 = -a_0 f_0.
$$
So $a_0 = b$.  Now (\ref{E:TeTf}) gives $a_n = b$, $n=0,1,2,\dots$.
So $T$ is a multiple of $T_1$ where 

\begin{equation}\label{E:HDS-example}
T_1 f_n = f_{n+1}, \qquad n=0,1,2,\dots.
\end{equation}

Suppose $T \in \A_1 \cap \A^\infty$.  Then $T f_0 = 0$ and
$T f_n = a_n f_{n-1}$ for some $a_n \in \CC$, $n=1,2, \dots$.  Since
$T_f \in \A_0$, by Lemma \ref{L:strong-Schur} there exists $b \in \CC$ such
that $T_f = -b \mathrm{I}$.  We compute
$$
T_f f_0 = -T R(f) f_0 = -\lambda a_1 f_0.
$$
So $a_1=b/\lambda$.  Now (\ref{E:TeTf}) gives
$$
a_n = \frac{nb}{\lambda+n-1}, \qquad n=1,2,\dots.
$$
So $T$ is a multiple of $T_1^*$.

Since $T_1$ and $T_1^*$ don't commute, it follows that at most one of
$\A_{-1}$ and $\A_1$ is nonzero.

Suppose $\A_{-1} \ne 0$.  Then $\A_1 = 0$ and so $\A_m=0$ for $m=1,2,\dots$
by Lemma \ref{L:down}.  Also $\A_{-1} = \CC T_1$ by the previous calculation
(and the fact that finite dimensional subspaces are closed).
Let $m>0$ and $T \in \A_{-m}$.  The equation $[T, T_1]=0$ implies
$T \in \CC T_1^m$.  By Lemma \ref{L:up}, $\A_{-m} \ne 0$, so in fact
$\A_{-m} = \CC T_1^m$.  So $\A$ is the strong-operator closure of the
algebra generated by $T_1$.

Suppose $\A_1 \ne 0$.  Then $\A_1 = 0$, and so $\A_m=0$ for $m=-1, -2, \dots$
by Lemma \ref{L:up}.  Also $\A_1 = \CC T_1^*$ by the previous calculation.
Let $m>0$ and $T \in \A_m$.  The equation $[T, T_1^*]=0$ implies
$T \in \CC (T_1^*)^m$.  By Lemma \ref{L:down}, $\A_m \ne 0$, so in fact
$\A_m = \CC (T_1^*)^m$.  So $\A$ is the strong-operator closure of the
algebra generated by $T_1^*$.

\smallskip
\noindent{\bf Case $I=\ZZ$:}

Suppose $T \in \A_{-1}$.  Then $T f_n = a_n f_{n+1}$ for some
$a_n \in \CC$, $n \in \ZZ$.  Since $T_e \in \A_0$, by
Lemma \ref{L:strong-Schur} there exists $b \in \CC$ such that
$T_e = -b \mathrm{I}$.  Together with (\ref{E:TeTf}), this
reads
$$
(\mu - n - 1) a_n - (\mu - n) a_{n-1} = -b,
$$
so there exists $a \in \CC$ such that
$$
a_n = \frac{a-bn}{\mu - n - 1}.
$$
Since $T$ commutes with $T_f$, we have
$$
\begin{aligned}
0
= & \; [T,T_f]f_n \\
= & \; (-(\lambda + \mu + n)a_{n+1} a_{n+2} +
      2 (\lambda + \mu + n + 1) a_n a_{n+2} -
      (\lambda + \mu + n + 2)a_n a_{n+1}) f_{n+3} \\
\therefore\quad 0
= & \; (a + b(1-\mu))(a + b (\lambda + \mu)).
\end{aligned}
$$
Therefore either $a=b(\mu-1)$ or $a=-b(\lambda+\mu)$.  So $T$ is a multiple
of either $T_2$ or $T_3$ where

\begin{equation}\label{E:PCS-example}
\begin{aligned}
T_2 f_n &= f_{n+1}, \qquad n \in \ZZ, \\
T_3 f_n &= \frac{\lambda+\mu+n}{n+1-\mu} f_{n+1}, \qquad n \in \ZZ.
\end{aligned}
\end{equation}

Suppose $T \in \A_1$.  Then $T f_n = a_n f_{n-1}$
for some $a_n \in \CC$, $n \in \ZZ$.  Since $T_f \in \A_0$, by
Lemma \ref{L:strong-Schur} there exists $b \in \CC$ such that
$T_f = -b \mathrm{I}$.  Together with (\ref{E:TeTf}), this
reads
$$
(\lambda + \mu + n - 1) a_n - (\lambda + \mu + n) a_{n+1} = -b,
$$
so there exists $a \in \CC$ such that
$$
a_n = \frac{a-bn}{\lambda + \mu + n - 1}.
$$
Since $T$ commutes with $T_e$, we have
$$
\begin{aligned}
0
= & \; [T,T_e]f_n \\
= & \; (-(\mu - n) a_{n-1} a_{n-2} + 
       2 (\mu - n + 1) a_n a_{n-2} -
       (\mu - n + 2) a_n a_{n-1})
       f_{n-3} \\
\therefore \quad 0
= & \; (a + b (\lambda + \mu - 1)) (a - b \mu).
\end{aligned}
$$
Therefore either $a=b (1 - \lambda - \mu)$ or $a=b \mu$.  So $T$
is a multiple of either $T_2^{-1}$ or $T_3^{-1}$.

Suppose $\A_{-1} \ne 0$.  Then by the previous calculation, $\A_{-1}=\CC T_2$
or $\A_{-1}=\CC T_3$.  Suppose $\A_{-1}=\CC T_2$.  Let $m>0$ and
$T \in \A_{-m}$.  The equation $[T,T_2]=0$ implies $T \in \CC T_2^m$.  By
Lemma \ref{L:up}, $\A_{-m} \ne 0$, so in fact $\A_{-m} = \CC T_2^m$.
Likewise, $\A_{-1}=\CC T_3$ implies $\A_{-m} =\CC T_3^m$ for all $m > 0$.

A similar argument shows that if $\A_1 \ne 0$ then either
$\A_m = \CC T_2^{-m}$ for all $m>0$ or
$\A_m = \CC T_3^{-m}$ for all $m>0$.

Since $T_2$ and $T_3$ don't commute unless they are equal (which occurs
precisely when $\mu = (1-\lambda)/2$), it now follows that $\A$ is the
strong operator closure of one of the following algebras:
$\CC [T_2]$, $\CC [T_3]$, $\CC [T_2^{-1}]$, $\CC [T_3^{-1}]$,
$\CC [T_2, T_2^{-1}]$ and $\CC [T_3, T_3^{-1}]$.

\smallskip
\noindent{\bf Case $R=D_\lambda^-$:}
Now $\A$ is $D_\lambda^-$-inductive iff it is $D_\lambda^+$-inductive.
However,
$\A_m$ and $\A_{-m}$ are interchanged.  To see this, let $(\kappa^+, \A^+)$
and $(\kappa^-, \A^-)$ denote the conjugation action of $\Mob$ on $\A$
under $D_\lambda^+$ and $D_\lambda^-$ respectively.  Then
$$
\begin{aligned}
\A^-_m
= & \; \{ T \in \A \st \kappa^-(g)T = \chi_m(g) T, \, g \in K \} \\
= & \; \{ T \in \A \st \kappa^+(g^*)T = \chi_m(g) T, \, g \in K \} \\
= & \; \{ T \in \A \st \kappa^+(g)T = \chi_m(g^*) T, \, g \in K \} \\
= & \; \{ T \in \A \st \kappa^+(g)T = \chi_{-m}(g) T, \, g \in K \} \\
= & \; \A^+_{-m}
\end{aligned}
$$
In particular, $\A^-_{-1} = \A^+_1 = \CC T_1^*$, where (we recall)

\begin{equation}\label{E:AHDS-example}
T_1^* f_n = \frac{n}{\lambda+n-1} f_{n-1}, \qquad n=1,2,\dots.
\end{equation}

\section{The homogeneous shifts}\label{S:HS}

Let $T \in \BH$ be homogeneous.  Recall (Definition 2.1 in \cite{Bagchi-Misra})
that a representation $R$ of $G$ on $\H$ is {\em associated} with $T$ if

\begin{equation}\label{E:associated}
\phi_g(T) = R(g)^{-1} T R(g), \qquad g \in G.
\end{equation}

Let $\A$ be the strong-operator closure of the subalgebra of $\BH$ generated
by $\{\phi(T) \st \phi \in \Mob \}$.

\begin{prop}\label{P:homogeneous-inductive}
The algebra $\A$ is $R$-inductive and $T \in \A_{-1}$.
\end{prop}

\begin{proof}
Since all elements of $\A$ are functions of a single operator, it is
abelian.  Furthermore, (\ref{E:associated}) shows that $\A$ is normalized
by $R(G)$.  So $\A$ is an $R$-inductive algebra.  If $\tilde\kappa$
and $\kappa$ are the conjugation representations of $G$ and $\Mob$ on
$\A$, then
$$
\begin{aligned}
\tilde\kappa(g)T &= R(g) T R(g)^{-1} = \phi_{g^{-1}}(T), & \qquad g \in G, \\
\therefore\quad \kappa(\phi)T &= \phi^{-1}(T), & \qquad \phi \in \Mob, \\
\therefore\quad \kappa(\phi_{\alpha,0})T &=
                \phi_{\alpha,0}^{-1}(T) = \alpha^{-1} T, &
                \qquad \alpha \in \TT,
\end{aligned}
$$
so $T \in \A_{-1}$.
\end{proof}

Let $I \in \{ \ZZ, \ZZ^+, \ZZ^- \}$.  Recall that an operator $T$ on a
Hilbert space $\H$ is called a {\em weighted shift} with weight sequence
$w_n$, $n \in I$, if there is a distinguished orthonormal basis
$x_n$, $n \in I$ such that $T x_n = w_n x_{n+1}$ for all $n \in I$.
The weighted shift is called a {\em bilateral shift}, {\em forward shift}
or {\em backward shift} according as $I=\ZZ$, $\ZZ^+$ or $\ZZ^-$.

It was shown in \cite{Bagchi-Misra} that if $T \in \BH$ is a homogeneous
weighted shift, then there is a unitary representation $R$ of
$G$ on $\H$ associated to $T$ and the $h$-isotypic subspaces of $\H$ are
precisely $\CC x_n$, i.e., for each $n \in I$ there exists
$\lambda_n \in \RR$ such that
$$
\{ F \in \H \st R(\exp(th))F = e^{-\lambda_n t}F, \, t \in \RR \} = \CC x_n.
$$
Further, it was shown that either $R$ is irreducible, or has the
decomposition $R \cong D_{2-\lambda}^- \oplus D_\lambda^+$.

Here, we will give different proofs of Lemma 5.1 and Theorem 5.2 from
\cite{Bagchi-Misra}.

First observe that $T_1$, $T_1^*$, $T_2$, and $T_3$ are homogeneous
by Theorem 2.3 in \cite{Bagchi-Misra}.

\begin{rem}\label{R:rigidity}
If both $T$ and $c T$ are associated to the same representation, then
$\phi(c T) = c \phi(T)$ for all $\phi \in \Mob$ and so $c=1$.
\end{rem}

Suppose $T$ is a homogeneous (weighted) shift and $R$ is associated to $T$.

Assume first that $R$ is irreducible.

If $R=D_\lambda^+$ is in the holomorphic discrete series, then $T$ is a
multiple of $T_1$ by (\ref{E:HDS-example}).  So, by Remark \ref{R:rigidity},
$T=T_1$.  If
$$
x_n = \sqrt{\frac{\Gamma(\lambda+n)}{\Gamma(1+n)}} f_n, \qquad n \in \ZZ^+,
$$
then $x_n$ is an orthonormal basis of $\H^{\lambda,0}$ and in terms of this
basis,
$$
\begin{aligned}
T_1 x_n
= & \; \sqrt{\frac{\Gamma(\lambda+n)}{\Gamma(1+n)}} T_1 f_n \\
= & \; \sqrt{\frac{\Gamma(\lambda+n)}{\Gamma(1+n)}} f_{n+1} \\
= & \; \sqrt{\frac{\Gamma(\lambda+n)}{\Gamma(1+n)}}
       \sqrt{\frac{\Gamma(2+n)}{\Gamma(\lambda+n+1)}}
       x_{n+1} \\
= & \; \sqrt{\frac{1+n}{\lambda+n}} x_n,
\end{aligned}
$$
so the weight sequence is $w_n = \sqrt{\frac{1+n}{\lambda+n}}$, $n \in \ZZ^+$.

If $R=D_\lambda^-$ is in the anti-holomorphic discrete series, then $T=T_1^*$
by (\ref{E:AHDS-example}) and Remark \ref{R:rigidity}.  If
$$
x_n = \sqrt{\frac{\Gamma(\lambda-n)}{\Gamma(1-n)}} f_{-n}, \qquad n \in \ZZ^-,
$$
then $x_n$ is an orthonormal basis of $\H^{\lambda,0}$ and in terms of this
basis,
$$
\begin{aligned}
T_1^* x_n
= & \; \sqrt{\frac{\Gamma(\lambda-n)}{\Gamma(1-n)}} T_1^* f_{-n} \\
= & \; \sqrt{\frac{\Gamma(\lambda-n)}{\Gamma(1-n)}}
       \frac{-n}{\lambda-n-1} f_{-n-1} \\
= & \; \sqrt{\frac{\Gamma(\lambda-n)}{\Gamma(1-n)}}
       \sqrt{\frac{\Gamma(-n)}{\Gamma(\lambda-n-1)}}
       \frac{-n}{\lambda-n-1} x_{n+1} \\
= & \; \sqrt{\frac{n}{1-\lambda+n}} x_{n+1}
\end{aligned}
$$
so the weight sequence is $w_n = \sqrt{\frac{n}{1-\lambda+n}}$, $n \in \ZZ^-$.

If $R=R_{\lambda,\mu}$ is in the principal series, then $T \in \{T_2, T_3\}$
by (\ref{E:PCS-example}) and Remark \ref{R:rigidity}.  Also $x_n = f_n$,
$n \in \ZZ$ is already an orthonormal basis of $\H^{\lambda,\mu}$.  So the
weight sequence is either $w_n = 1$, $n \in \ZZ$, or
$w_n = \frac{\lambda+\mu+n}{n+1-\mu}$, $n \in \ZZ$.

If $R=R_{\lambda,\mu}$ is in the complementary series, then
$T \in \{T_2, T_3\}$ by (\ref{E:PCS-example}) and Remark \ref{R:rigidity}.
If
$$
x_n = \sqrt{\frac{\Gamma(\lambda+\mu+n)}{\Gamma(1-\mu+n)}} f_n,
\qquad n \in \ZZ,
$$
then $x_n$ is an orthonormal basis of $\H^{\lambda,0}$ and in terms of this
basis,
$$
\begin{aligned}
T_2 x_n
= & \; \sqrt{\frac{\Gamma(\lambda+\mu+n)}{\Gamma(1-\mu+n)}} T_2 f_n \\
= & \; \sqrt{\frac{\Gamma(\lambda+\mu+n)}{\Gamma(1-\mu+n)}} f_{n+1} \\
= & \; \sqrt{\frac{\Gamma(\lambda+\mu+n)}{\Gamma(1-\mu+n)}}
       \sqrt{\frac{\Gamma(2-\mu+n)}{\Gamma(\lambda+\mu+n+1)}}
       x_{n+1} \\
= & \; \sqrt{\frac{1-\mu+n}{\lambda+\mu+n}} x_n,
\end{aligned}
$$
so the weight sequence is $w_n = \sqrt{\frac{1-\mu+n}{\lambda+\mu+n}}$,
$n \in \ZZ$.

Now suppose $R = D_{2-\lambda}^- \oplus D_\lambda^+$ acting on
$\H = \H^{(2-\lambda),0} \oplus \H^{\lambda,0} = \H_1 \oplus \H_2$ (say).  Let
$$
g_n =
\begin{cases}
 (f_{-1-n},0), & n < 0, \\
 (0, f_n), & n \ge 0.
\end{cases}
$$
Then $g_n$, $n \in \ZZ$ is a basis of $\H$,
$$
\{ F \in \H \st R(\exp(th))F = e^{-i(2n+\lambda)t} F \} = \CC g_n,
\qquad n \in \ZZ,
$$
and $T g_n = a_n g_{n+1}$ for some $a_n \in \CC$, $n \in \ZZ$.
For future reference, we record

\begin{equation}\label{E:ReRf}
\begin{aligned}
R(e) g_n & =
\begin{cases}
(1-\lambda-n) g_{n-1}, & n < 0, \\
0, & n=0, \\
-n g_{n-1}, & n>0,
\end{cases} \\
R(f) g_n & =
\begin{cases}
(n+1) g_{n+1}, & n<-1, \\
0, & n=-1, \\
(\lambda+n) g_{n+1}, & n > -1.
\end{cases}
\end{aligned}
\end{equation}

Write
$$
T=\(
\begin{matrix}
T_{11} & T_{12} \\
T_{21} & T_{22}
\end{matrix}
\)
$$
where $T_{ij}:\H_j \to \H_i$, $i,j \in \{ 1,2 \}$, $T_{12}=0$ and
$$
T_{21} f_n =
\begin{cases}
0 & n>0 \\
r f_0 & n=0,
\end{cases}
$$
for some $r \in \CC$.  So $D_{2-\lambda}^-$ is associated to $T_{11}$ and
$D_\lambda^+$ is associated to $T_{22}$, and so
$$
\begin{aligned}
T_{11} f_0 &= 0, \\
T_{11} f_n &= \frac{n}{1-\lambda+n} f_{n-1}, \qquad n=1,2,\dots, \\
T_{22} f_n &= f_{n+1}, \qquad n=0,1,2,\dots.
\end{aligned}
$$
So
$$
a_n =
\begin{cases}
\frac{(1 + n)}{\lambda+n}, & n<-1, \\
r, & n=-1, \\
1, & n>-1.
\end{cases}
$$

We claim that $\lambda=1$.  This requires the full strength of the condition
that $T$ is homogeneous (at the infinitesimal level).  We compute
$$
\begin{aligned}
\kappa(L) T
= & \; \left. \frac{d}{dt} \right|_{t=0} \kappa(\exp tL) T \\
= & \; \left. \frac{d}{dt} \right|_{t=0} \phi_{1, -\tanh t}^{-1}(T) \\
= & \; \left. \frac{d}{dt} \right|_{t=0} \phi_{1, \tanh t}(T) \\
= & \; \left. \frac{d}{dt} \right|_{t=0}
       \frac{T - (\tanh t) \mathrm{I}}{I - (\tanh t) T} \\
= & \; T^2 - \mathrm{I},
\end{aligned}
$$
and
$$
\begin{aligned}
\kappa(M) T
= & \; \left. \frac{d}{dt} \right|_{t=0} \kappa(\exp tM) T \\
= & \; \left. \frac{d}{dt} \right|_{t=0} \phi_{1, -i\tanh t}^{-1}(T) \\
= & \; \left. \frac{d}{dt} \right|_{t=0} \phi_{1, i\tanh t}(T) \\
= & \; \left. \frac{d}{dt} \right|_{t=0}
       \frac{T - i(\tanh t) \mathrm{I}}{I + (\tanh t) T} \\
= & \; -i (T^2 + \mathrm{I}).
\end{aligned}
$$
Therefore $\kappa(e) T = - \mathrm{I}$ and
$\kappa(f) T = T^2$.  Now, from (\ref{E:ReRf}) we see that
$(\kappa(f)T) g_{-1} = [R(f), T] g_{-1} = r \lambda g_1$.  But
$T^2 g_{-1} = r g_1$.  So $\lambda=1$, proving the claim.

In this situation, $x_n=g_n$, $n \in \ZZ$ is already an orthonormal basis
of $\H$.  So the weight sequence is
$$
w_n =
\begin{cases}
1, & n<-1, \\
r, & n=-1, \\
1, & n>-1.
\end{cases}
$$

\bibliographystyle{amsplain}
\bibliography{v15-iahs}

\end{document}